\documentclass[11pt]{article}
\usepackage{amsthm}
\usepackage{latexsym, amssymb, amsmath}
\usepackage{graphicx}
\usepackage[dvipsnames]{xcolor}

\begin{document}
\numberwithin{equation}{section}

\newtheorem{THEOREM}{Theorem}
\newtheorem{PRO}[THEOREM]{Proposition}
\newtheorem{XXXX}{\underline{Theorem}}
\newtheorem{CLAIM}{Claim}
\newtheorem{COR}{Corollary}
\newtheorem{LEMMA}{Lemma}
\newtheorem{REM}{Remark}
\newtheorem{EX}{Example}
\newenvironment{PROOF}{{\bf Proof}.}{{\ \vrule height7pt width4pt depth1pt} \par \vspace{2ex} }
\newcommand{\Bibitem}[1]{\bibitem{#1} \ifnum\thelabelflag=1 
  \marginpar{\vspace{0.6\baselineskip}\hspace{-1.08\textwidth}\fbox{\rm#1}}
  \fi}
\newcounter{labelflag} \setcounter{labelflag}{0}
\newcommand{\labelon}{\setcounter{labelflag}{1}}
\newcommand{\Label}[1]{\label{#1} \ifnum\thelabelflag=1 
  \ifmmode  \makebox[0in][l]{\qquad\fbox{\rm#1}}
  \else\marginpar{\vspace{0.7\baselineskip}\hspace{-1.15\textwidth}\fbox{\rm#1}}
  \fi \fi}

\newcommand{\LEFTLINE}{\ifhmode\newline\else\noindent\fi}
\newcommand{\RIGHTLINE}[1]{\LEFTLINE\rightline{#1}}
\newcommand{\CENTERLINE}[1]{\LEFTLINE\centerline{#1}}
\def\BOX #1 #2 {\framebox[#1in]{\parbox{#1in}{\vspace{#2in}}}}
\parskip=8pt plus 2pt
\def\AUTHOR#1{\author{#1} \maketitle}
\def\Title#1{\begin{center}  \Large\bf #1 \end{center}  \vskip 1ex }
\def\Author#1{\vspace*{-2ex}\begin{center} #1 \end{center}  
 \vskip 2ex \par}
\renewcommand{\theequation}{\arabic{section}.\arabic{equation}}
\def\bdk#1{\makebox[0pt][l]{#1}\hspace*{0.03ex}\makebox[0pt][l]{#1}\hspace*{0.03ex}\makebox[0pt][l]{#1}\hspace*{0.03ex}\makebox[0pt][l]{#1}\mbox{#1} }
\def\psbx#1 #2 {\mbox{\psfig{file=#1,height=#2}}}


\newcounter{firstbib}
\newcommand{\FG}[2]{{\includegraphics[height=#1mm]{#2.eps}}}

 
 
\Title{A New Family of Nonnegative Sine Polynomials}

\begin{center}
MAN KAM KWONG
\end{center}

\begin{center}
\emph{Department of Applied Mathematics\\ The Hong Kong Polytechnic University,\\ Hunghom, Hong Kong}\\
\tt{mankwong@polyu.edu.hk}\\[4ex]
\end{center}

\par\vspace*{\baselineskip}\par

\newcommand{\mb}{\mathbf}
\newcommand{\Cr}{\color{red}}
\newcommand{\tx}[1]{\mbox{ #1 }}
\newcommand{\dc}{\searrow}
\newcommand{\bR}{\vspace*{5pt}\begin{REM}\em}
\newcommand{\eR}{\qed\end{REM}\vspace*{5pt}}
\newcommand{\irp}[1]{\par \hfill IRTP: \fbox{ \em #1 }}

\begin{abstract}
\parskip=6pt
We study sine polynomials of the form
$$  \kappa \,\sin(x) + \sum_{k=2}^{n-1} \sin(kx)  + \lambda \,\sin(nx) \qquad  (\kappa , \lambda  \in \mathbb R )  $$
that is nonnegative for all $ x\in[0,\pi ] $. These include, in particular,
$$  \frac{5}{4} \,\sin(x) + \sum_{k=2}^{n-1} \sin(kx)  + \frac{2n-3}{4n} \,\sin(nx)  )  \qquad  (n \tx{odd}),  $$
and
$$  \left( n + \frac{1}{2} \right)  \,\sin(x) + \sum_{k=2}^{n} \sin(kx)  \hspace*{24mm} (n \tx{even}).  $$

We also characterize all nonnegative sine polynomials of degree 3 and all nonnegative cosine
polynomials of degree 2.
\end{abstract}



\par\vspace*{2\baselineskip}\par

{\bf{Acknowledgments.}} This research was inspired while the author
was visiting the National Center for Theoretical Sciences
(No.\ 1 Sec.\ 4 Roosevelt Rd.,
National Taiwan University,
Taipei, 106, Taiwan) in January, 2016. The author is grateful for
the support and hospitality of the Center.

The research of this author is partially supported by the Hong Kong Government GRF Grant PolyU 5012/10P and the Hong Kong Polytechnic University Grants G-YK49 and G-U751.

\newpage
\section{Introduction}

Nonnegative trigonometric polynomials are useful. Many excellent surveys on their history
and applications are available in the 
literature. See, for example, Brown \cite{b1}, 
Dumitrescu \cite{d1}, and Koumandos \cite{ko1}, and the many references therein.

\noindent
{\bf Notations}. For convenience, we adopt the following notations:
$$  c_k = c(k) := \cos(kx) \, ,  $$
$$  s_k = s(k) := \sin(kx) \, ,  $$
$$  [ a_1,a_2,...,a_n]_s := \sum_{k=1}^{n} a_ks_k \,,  $$
$$  [ a_0,a_1,...,a_n]_c := \sum_{k=0}^{n} a_kc_k \,,  $$
where $ \left\{ a_k\right\} _{1}^n $ is a given sequence of real numbers. We use the acronym NN
as an abbreviation for either the adjective ``nonnegative'' or the ``property of
being nonnegative'', depending on the context. The property usually concerns
an expression depending on some variable $ x $ in some interval $ I $.
When the interval $ I $ is not explicitly specified, the default is $ [0,\pi ] $.
We also abbrviate ``lefthand (righthand) side'' to RHS (LHS).

In 1958, Vietoris \cite{v} proved a deep result concerning NN trigonometric polynomials. It was
improved by Belov \cite{be} in 1995.

\begin{THEOREM}[Vietoris-Belov]
Suppose that $ \left\{ a_k\right\} _0^\infty  $ is a sequence of non-increasing positive numbers satisfying 
the condition 
$$  \sum_{k=1}^{m} (-1)^{k-1} k\,a_k \geq  0 , \quad  \tx{for all (even)} m \geq  1 .  \eqno (\rm{B})  $$
Then for all $ x\in[0,\pi ] $ and $ m\geq 1 $,
$$  \sum_{k=0}^m a_k\,c_k \geq 0  \qquad  \tx{and} \qquad   \sum_{k=1}^m a_k\,s_k \geq 0 .   $$

For sine polynomials with $ a_k\searrow $, {\rm (B)} is also necessary for all its
partial sums to be NN.
\end{THEOREM}

\bR
Frequently, Theorem 1 is applied to a finite sequence $ \left\{ a_k\,:\,k\leq n\right\}  $
(simply take $ a_k=0 $ for $ k>n $). The length of the sequence $ n $ is called
the degree of the trigonometric polynomial.
Then (B) consists of $ \frac{n}{2} $ (if
$ n $ is even) or $ \frac{n-1}{2} $ inequalities.
\eR

\bR
The original Vietoris result covers a smaller family of polynomials including

$ \hspace*{20mm}   a_0 = a_1 = 1, \qquad     \textstyle  a_{2k+1} = a_{2k} = \frac{1\cdot3\cdot5\cdots (2k-1)}{2\cdot4\cdot6\cdots (2k)} \quad  (k=1,2,...) . $
\eR

\bR
Vietoris' (but not Belov's) result has recently been extended by the author
\cite{kw3} to sine polynomials with non-decreasing and non-decaying coefficients.
\eR

Contrast Theorem~1 with the following NN criterion, 
due to Fej\'er \cite{f1}, in which
only the NN of the entire polynomial, not its partial sums, is 
asserted.

\begin{THEOREM}[Fej\'er]
Suppose $ a_k\dc $, $ {}>0 $, $ k\leq n $, satisfy
$$  a_{k}+a_{k+2}\geq 2a_{k+1}, \qquad  k=1,2,...,n-2.  $$
Then
$$  \sum_{k=1}^{n-1} a_k\,s_k + \frac{a_n}{2} \,s_n \quad  \tx{is NN in} [0,\pi ].  $$
\end{THEOREM}

\bR
Among the simplest examples are
$ [n,n-1,n-2, ..., 1]_s $
(attributed to E. Luk\'acs by Fej\'er \cite{f1}) and
\begin{equation}  \sigma (x) = [\,1,1, ...,1, \textstyle  \frac{1}{2} \,]_s .  \Label{s21}  \end{equation}
None of their proper partial sums are NN.

New NN polynomials can be constructed by taking positive linear combinations of
known ones. For instance,
$$  [2,1]_s + [2,2,2,2,1]_s = [4,3,2,2,1]_s  $$
is NN; but it satisfies neither Theorem 1 nor Theorem 2.
\eR
  
\bR
Strangely enough, no cosine analog of Theorem~2 is known.
\eR


A long term goal of the author is to establish general NN criteria that
extend or unify these known results. In the meantime, the discovery of
new families of NN polynomials will help to achieve that goal. It is the
purpose of this article to present such a new family that includes
$ \sigma  $ as a particular case.

In Section 2, we show that, for odd $ n $,
\begin{equation}  \phi (x) = \textstyle  [\frac{5}{4} ,1,1,\cdots ,1,\frac{2n-3}{4} ]_s,   \Label{phi}  \end{equation}
is NN. Here, all the middle coefficients are 1.
The proof is non-trivial. In Section 3, we consider
polynomials of the more general form $ [\kappa ,1,1,\cdots ,1,\lambda ]_s $.

\bR
Belov's criterion (B) imposes multiple inequalities on the coefficients and
reaps multiple (all partial sums) NN results. It is natural to wonder
whether imposing only a single inequality, namely, the $ n $-th one 
(i.e.\ (\ref{Bn3}) below),
can lead to any useful conclusion. It turns out that (\ref{Bn3}) alone
is not enough to guarantee NN.
However, it is known that (\ref{Bn3}) is necessary for NN.
In fact, the following more general assertion holds.
\eR

\newpage

\begin{PRO} Let $ \left\{ a_k\right\} _1^n $ be a finite sequence of real numbers, positive or negative.
\begin{enumerate}
\item[\rm{(i)}] A necessary condition for 
$ [a_1,...,a_n]_s $
to be NN in a right neighborhood of $ x=0 $ is
\begin{equation}  \sum_{k=1}^{n}  ka_k \geq  0 .  \Label{Bn1}  \end{equation}
In case it is known that $ \sum ka_k=0 $,
an additional necessary condition is
\begin{equation}  \sum_{k=1}^{n}  k^3a_k \leq  0 .  \Label{Bn2}  \end{equation}
\item[\rm{(ii)}] A necessary condition for 
$ [a_1,...,a_n]_s $
to be NN in a left neighborhood of $ x=\pi  $ is
\begin{equation}  \sum_{k=1}^{n}  \, (-1)^{k+1}ka_k \geq  0 .  \Label{Bn3}  \end{equation}
In case it is known that 
\begin{equation}  \sum \, (-1)^{k+1}ka_k=0,  \Label{Bn0}  \end{equation}
an additional necessary condition is
\begin{equation}  \sum_{k=1}^{n} \, (-1)^{k+1} k^3a_k \leq  0 .  \Label{Bn4}  \end{equation}
\end{enumerate}
\end{PRO}

\noindent
{\Cr [In an earlier version, (\ref{Bn2}) and (\ref{Bn0}) had the typo where $ \leq  $
was mistakenly typed as $ \geq  $.]}

\begin{PROOF}
We only give the proof of (i), that of (ii) being similar. Denote 
$ f(x) = \sum a_ks_k.  $
Since $ f(0)=0 $ and $ f(x)\geq 0 $ for $ x\in(0,\epsilon ) $ for some $ \epsilon >0 $, we have
$ f'(0)\geq 0 $. Note that $ f'(x)=\sum ka_kc_k $. Substituting $ x=0 $ in this
identity gives (\ref{Bn1}).

If it happens that $ \sum ka_k=0 $, then
$ f'(0)=0 $. Note that $ f''(x)=-\sum k^2a_ks_k $ and 
$ f'''(x)=-\sum k^3a_kc_k $. Hence, $ f''(0)=0 $ and $ f'''(0)=-\sum k^3a_k $.
The only way that $ f(x) $ can be NN in $ (0,\epsilon ) $ is to have $ f'''(0)\geq 0 $,
yielding the desired necessary condition (\ref{Bn2}).
\end{PROOF}

\bR
$ \sigma  $ of (\ref{s21}) satisfies (\ref{Bn0}) and (\ref{Bn4}) when $ n $ is even. When $ n $ is odd, it
satisfies (\ref{Bn3}) with strict inequality. On the other hand, $ \phi  $ of (\ref{phi}) satisfies
(\ref{Bn0}) and (\ref{Bn4}) (with equality) when $ n $ is odd (see (\ref{k3}) below).
\eR

\bR
It is easy to show that, after replacing $ \geq  $ by $ > $,
condition (\ref{Bn1}) (or (\ref{Bn3})) is sufficient for the sine polynomial
to be NN in some right (left) neighborhood of $ x=0 $ ($ x=\pi  $). 
However, this simple fact is not
too useful because there is no information on how large this neighborhood can be.
\eR

We will make use of the well-known inequalities (for $ t\geq 0 $) 
\begin{equation}  \hspace*{20mm}      \sin(t)\leq t, \qquad  \cos(t) \geq  1 - \frac{t^2}{2}  \Label{sc0}  \end{equation}
\begin{equation}  \hspace*{22mm}       t - \frac{t^3}{6} \,\, \leq  \,\, \sin(t) \,\, \leq  \,\, t - \frac{t^3}{6} + \frac{t^5}{120}  \Label{sin1}  \end{equation}
\begin{equation}  1 - \frac{t^2}{2} + \frac{t^4}{24} - \frac{t^6}{720} \,\, \leq  \,\, \cos(t) \,\, \leq  \,\, 1 - \frac{t^2}{2} + \frac{t^4}{24}   \Label{cos1}  \end{equation}
in the next section. These can be proved by iteratively
integrating the simple inequality $ \cos(t)\leq 1 $.
The polynomials involved are truncated Taylor series
of $ \sin $ and $ \cos $, respectively.

The Sturm procedure is a useful technique for proving the NN of specific 
trigonometric polynomials with given constant coefficients. In a nutshell, the NN of a given trigonometric polynomial is equivalent
to that of a corresponding algebraic polynomial. The latter can be verified using
the classical Sturm result on the exact number of real solutions of an algebraic polynomial
in a given interval of real numbers, with the help of a computer software.
In our study, we used the extremely helpful mathematical software 
MAPLE 2016. Detailed expositions can be found in 
\cite{kw3}, \cite{ak1}, and \cite{kw4}.

\section{$\phi  = [\,\frac{5}{4} ,1,1, ... , 1, \frac{2n-3}{4n} \,]_s \,$, $ n\geq 3 $, Odd }

Other than the first and last, all coefficients of $ \phi  $ are 1.

\begin{LEMMA}
For odd $ n\geq 3 $, $\phi $ is NN in $ [0,\pi ] $.
\end{LEMMA}  \begin{PROOF}
The Sturm procedure have been used to confirm the conclusion
for $ n=3,5,7 $, and $ 9 $.
In the following, we assume that $ n\geq 11 $.

In $ \big[0,\frac{\pi }{n} \big] $, all terms of the polynomial in question are NN. Hence, the sum 
  is positive. It remains to establish NN in the remaining interval $ \big[\frac{\pi }{n} ,\pi \big] $.
  By applying the reflection $ x\mapsto\pi -x $, we see that this
  is equivalent to the NN, in $ \big[0,\pi -\frac{\pi }{n} \big] $, of
$$  \textstyle    S_1 = \left[  \, \frac{5}{4} ,-1,1,-1,1, ... \, , -1,\frac{2n-3}{4n} \, \right] _s  .  $$

The product-to-sum identity gives
\begin{eqnarray}
   8c{\textstyle  \left( \frac{1}{2} \right) } S_1 &=&  5  s{\textstyle  \left( \frac{1}{2} \right) } +  s{\textstyle  \left( \frac{3}{2} \right) } - \frac{2n+3}{n} \, s{\textstyle  \left( n - \frac{1}{2} \right) } + \frac{2n-3}{n} \, s{\textstyle  \left( n + \frac{1}{2} \right) } \nonumber \\
    &=&  5  s{\textstyle  \left( \frac{1}{2} \right) } +  s{\textstyle  \left( \frac{3}{2} \right) } + 4 s{\textstyle  \left( \frac{1}{2} \right) } c(n) - \frac{6}{n} \,  c{\textstyle  \left( \frac{1}{2} \right) } s(n) .  \Label{s1}
\end{eqnarray}
By substituting $ 2x $ for $ x $ in (\ref{s1}),
we see that the desired assertion is the same as the NN of
\begin{equation}  S_2 =  5  s_1 + s_3 + 4s_1c_{2n} - \frac{6}{n} \, c_1s_{2n}     \Label{s2}  \end{equation}
in $ \big[0,\frac{\pi }{2} -\frac{\pi }{2n} \big] $.
We establish this claim separately in each of five subintervals:
\begin{equation}  \bigcup_{k=1}^5 I_k =  [\,0.14,1.4\,] \cup \big[1.4,{\textstyle  \frac{\pi }{2} - \frac{\pi }{2n} \big] \cup \big[\frac{\pi }{2n} ,0.14\big] \cup \big[ \frac{\pi }{4n} ,\frac{\pi }{2n} \big] \cup \big[0 ,\frac{\pi }{4n} \big] . }  \end{equation}
The following numbered paragraphs correspond to these subintervals in 
the same order.

\begin{enumerate}
\item Since $ c_{2n}\leq 1 $ and $ s_{2n}\leq 1 $, (\ref{s2}) yields
\begin{equation}  S_2 \geq  S_3 := s_1 + s_3 - \frac{6}{n} \, c_1 ,  \Label{s3}  \end{equation}
for all $ n $. In particular, with the choice $ n=11 $,
\begin{equation}  S_2 \geq  s_1 + s_3 - \frac{6}{11} \, c_1 .  \Label{s4}  \end{equation}
Using the Sturm procedure for general trigonometric polynomials, as described in \cite{kw4},
we can verify that the RHS of (\ref{s4}) is NN
in $ I_1 $.

\item Differentiating (\ref{s3}) gives
\begin{equation}  S_3'(x) = c_1+3\,c_3 +\frac{6}{n} \,s_1.   \Label{s5}  \end{equation}
At $ x=\frac{\pi }{2} -\frac{\pi }{2n}  $,
\begin{eqnarray}
         S_3' &=& \sin{\textstyle  \left(  \frac{\pi }{2n} \right) } - 3\,\sin{\textstyle  \left(  \frac{3\pi }{2n} \right) } + \frac{6}{n} \,\cos{\textstyle  \left(  \frac{\pi }{2n} \right) } \nonumber  \\
	 &\leq & \sin{\textstyle  \left(  \frac{\pi }{2n} \right) } - 3\,\sin{\textstyle  \left(  \frac{3\pi }{2n} \right) }+ \frac{6}{n} \,.  \Label{s6}
\end{eqnarray}
Using the first inequality of (\ref{sc0}) for the first term and the LHS of
(\ref{sin1}) for the second term,
we can show that the last expression of (\ref{s6}) is $ {}\leq 0 $, implying
that $ S_3'\big(\frac{\pi }{2} -\frac{\pi }{2n} \big)\leq 0 $.
Differentiating (\ref{s5}) gives
$$  S_3''(x) = (-s_1-9s_3) +\frac{6}{n} \, c_1.  $$
In $ I_2 $, the two terms in the RHS of the above expression are
nonnegative. Thus, $ S_3''(x)>0 $, implying that
$ S_3' $ is increasing in this interval. Consequently, $ S_3' $ is nonnegative in 
$ I_2 $, implying that $ S_3 $ is decreasing there.
At the right endpoint,
\begin{eqnarray}
S_3 &=& \cos{\textstyle  \left(  \frac{\pi }{2n} \right) } - \cos{\textstyle  \left(  \frac{3\pi }{2n} \right) } - \frac{6}{n} \,\sin{\textstyle  \left(  \frac{\pi }{2n} \right) } \nonumber  \\
    &\geq & \cos{\textstyle  \left(  \frac{\pi }{2n} \right) } - \cos{\textstyle  \left(  \frac{3\pi }{2n} \right) } - \frac{3\pi }{n^2} \,.    \Label{s7}
\end{eqnarray}
Using the appropriate inequality of (\ref{sc0}) for the first terms 
and of (\ref{cos1})
for the second term, we can show that the
expression in the second line  of (\ref{s7}) is NN. Hence,
$ S_3 $ is NN at $ x=\frac{\pi }{2} -\frac{\pi }{2n}  $. Since $ S_3 $ is decreasing,
it is NN in $ I_2 $. By (\ref{s3}), $ S_2 $ is NN in the
same interval. 

\item Using (\ref{s5}), we can easily verify
that $ S_3 $ is increasing in $ [0,0.14] $. At $ x=\frac{\pi }{2n} $,
\begin{eqnarray*}
S_3 &=& \sin{\textstyle  \left(  \frac{\pi }{2n} \right) } + \sin{\textstyle  \left(  \frac{3\pi }{2n} \right) } - \frac{6}{n} \,\cos{\textstyle  \left(  \frac{\pi }{2n} \right) } \\
&\geq & \sin{\textstyle  \left(  \frac{\pi }{2n} \right) } + \sin{\textstyle  \left(  \frac{3\pi }{2n} \right) } - \frac{6}{n} \,.  
\end{eqnarray*}
Using (\ref{sc0}), we can show that the last expression, and hence, also $ S_3 $
is NN.

\item We estimate the first two terms in (\ref{s2}) using (\ref{sin1}):
$$  5\,s_1 + s_3 \geq  8x - \frac{16}{3} \,x^3 .  $$
  In this interval, the third term in (\ref{s2}) is $ {}\leq 0 $, and 
$$  |c_{2n} | = |\cos(2nx)| = \sin\big(2nx - {\textstyle \frac{\pi }{2} } \big) \leq  2nx - {\textstyle \frac{\pi }{2} } .  $$
Hence
$$  4s_1c_{2n} \geq  -4x|\cos(2nx)| \geq  -4x \big( 2nx - {\textstyle \frac{\pi }{2} } \big) .  $$
To estimate the fourth term, note that $ c_1\leq 1 $ and
$$  s_{2n}  =  \cos(y) \leq  1 - \frac{y^2}{2} + \frac{y^4}{24} \,,  $$
where $ y=2nx-\frac{\pi }{2} \in[0,\frac{\pi }{2} ] $. Applying these estimates to (\ref{s2}), we obtain
$$  \textstyle  n S_2 \geq  - \frac{1}{4} \, y^4 -\frac{2}{3n^2} \,y^3 +  \left[  1 - \frac{\pi }{n^2} \right] \,y^2 + \left[  4 -\pi  - \frac{\pi ^2}{2n^2} \right]  \,y + \left[ 2\pi  -6 - \frac{\pi ^3}{12n^2} \right]  .  $$
For any fixed $ y\in\big[0,\frac{\pi }{2} \big] $, the RHS is increasing in $ n $.
The Sturm procedure can be applied to show that, for
the choice $ n=11 $, the RHS is NN in $ \big[0,\frac{\pi }{2} \big] $. Hence, $ S_2 $
is NN for all $ n\geq 11 $.

\item For this last subinterval, we apply the transformation: $ x\mapsto\alpha x $,
  $ \alpha =\frac{1}{2n} \in\big(0,\frac{1}{22} \big] $, to $ S_2 $.
  The desired assertion is equivalent to the NN of
$$  S_5 = 5 s(\alpha ) +  s(3\alpha ) + 4c(1) s(\alpha )  - 12\alpha  s(1)c(\alpha ).   $$
Applying the LHS of (\ref{sin1}) and (\ref{cos1}) to the first three terms, 
and the RHS of the inequalities to the last term, we obtain a lower bound
for 
$$  \frac{S_2}{\alpha x^5} \geq  \left[ {\frac {\alpha ^2}{1080}}-{\frac {\alpha ^4}{240}} \right] { x}^{4}+ \left[ \frac{\alpha ^2}{45} -{\frac {1}{180}}+ \frac{\alpha ^4}{12} \right] {x }^{2} + \left[  \frac{1}{15} - \frac{2\alpha ^2}{3} - \frac{\alpha ^4}{2} \right]  .  $$
By ignoring those positive terms involving $ \alpha  $ and estimating the negative ones
using $ \alpha \leq \frac{1}{22} $ and $ x\leq \frac{\pi }{2} $, we see that 
\begin{eqnarray*}
  \frac{S_2}{\alpha x^5}  &\geq &  -\,{\frac {\alpha ^4}{240}} \, { x}^{4} -\, {\frac {1}{180}} \, {x }^{2} + \left[  \frac{1}{15} - \frac{2\alpha ^2}{3} - \frac{\alpha ^4}{2} \right]  \\
    &\geq &  -\,{\frac {\big(\frac1{22}\big)^4}{240}} \, { \left( \frac{\pi }{2} \right) }^{4} -\, {\frac {1}{180}} \, {\left( \frac{\pi }{2} \right)  }^{2} +  \frac{1}{15} - \frac{2\big(\frac1{22}\big)^2}{3} - \frac{\big(\frac1{22}\big)^4}{2}  \quad  > \quad  0 .
\end{eqnarray*}
\end{enumerate}
The proof of the Lemma is now complete.
\end{PROOF}

\section{$[\kappa  ,\,1,\,1,\, ... , \,1,\, \lambda  ]_s \,$, $ n\geq 3 $}

In this section, we present our main result which concerns the family of sine 
polynomials of the form given in the section heading, with $ \kappa ,\lambda \in\mathbb R $,
and all other coefficients 1. Define
$$  {\cal P}_n = \left\{  (\kappa ,\lambda ) : [\kappa ,1,1,\cdots ,1,\lambda ]_s \tx{is NN for} x \in [0,\pi ] \right\} .  $$
When $ n $ is clear from the context, we suppress the subscript and write simply $ \cal P $.

\begin{THEOREM}
For any $ \lambda  $, there exists $ \kappa _0=\kappa _0(\lambda ;n) $ such that $ (\kappa ,\lambda )\in{\cal{P}}_n $ iff
$ \kappa \geq \kappa _0 $.
\begin{itemize}
\item[\rm(i)] For odd $ n $,
\begin{eqnarray}
        \kappa _0 &=& \textstyle  \frac{n+1}{2} -n\lambda  , \hspace*{9mm}  \lambda  \in \big(-\infty ,\frac{2n-3}{4n} \,\big] \Label{k01} \\
	\kappa _0 &>& \textstyle  \frac{n+1}{2} -n\lambda  , \hspace*{9mm}  \lambda  \in \big(\frac{2n-3}{4n} ,\frac{1}{2} \,\big] \Label{k02} \\
	\kappa _0 &>& 1                , \hspace*{22.8mm} \textstyle  \lambda  \in \big(\frac{1}{2} ,\infty \big) \, .  \Label{k03}
\end{eqnarray}
\item[\rm(ii)] For even $ n $,
\begin{eqnarray}
        \kappa _0 &>& \textstyle  1  , \hspace*{23mm}  \lambda  \in \big(-\infty ,\frac{1}{2} \big) \hspace*{6mm} \Label{k04} \\
	\kappa _0 &=& \textstyle  n\lambda  - \frac{n-1}{2} , \hspace*{9mm} \lambda  \in\big[\,\frac{1}{2} ,\infty \big) \, . \Label{k05}
\end{eqnarray}
\end{itemize}
\end{THEOREM}

\par\vspace*{\baselineskip}\par

\begin{figure}[h]
\begin{center}
\FG{50}{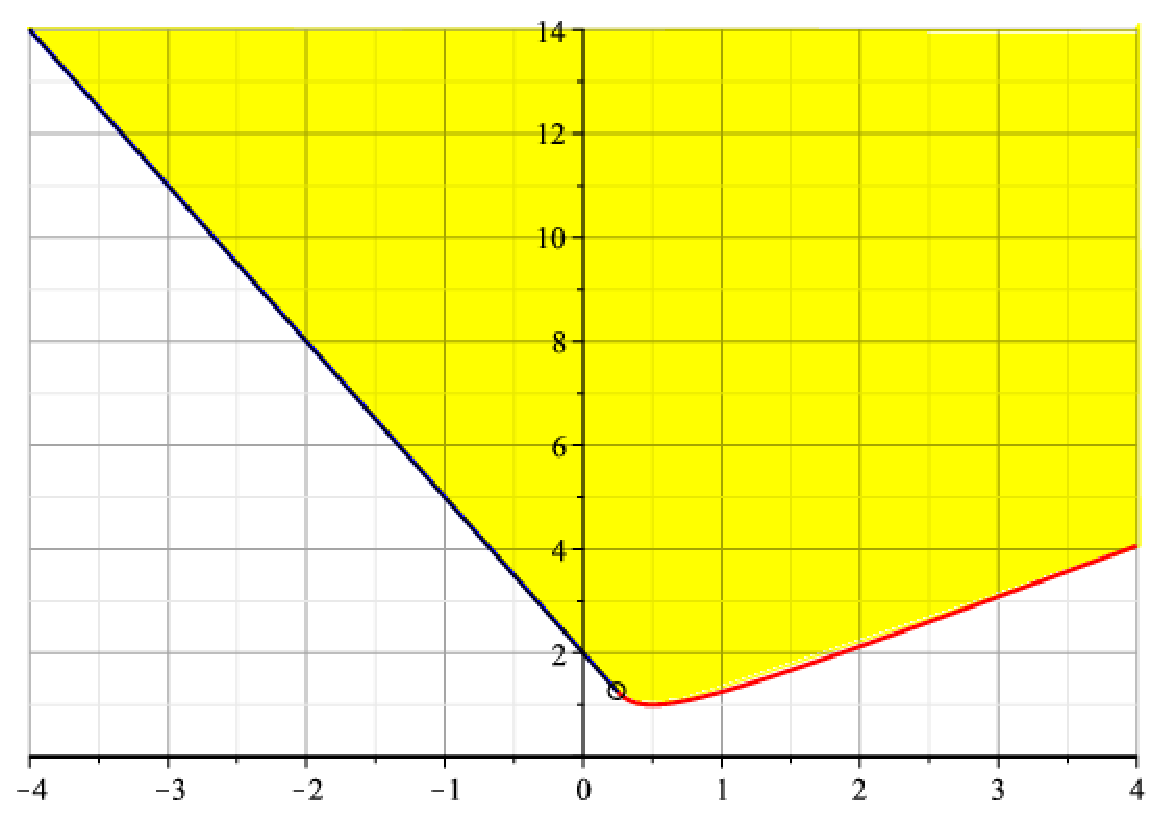}
\end{center}
\par\vspace*{-19mm}\par \hspace*{77mm} {\footnotesize  $ A $}
\par\vspace*{-0.85mm}\par   \hspace*{78.5mm} {\scriptsize  $ \circ $}
\par\vspace*{-2.9mm}\par   \hspace*{80.2mm} {\footnotesize  $ B $}
\par\vspace*{-26mm}\par \hspace*{90mm} $ {\cal{P}}_3 $
\par\vspace*{32mm}\par
\caption{${\cal{P}}_3$, $ A=\big(\frac{2n-3}{4n} ,\frac{5}{4} \big) $, $ B=\big(\frac{1}{2} ,1\big) $}
\par\vspace*{15mm}\par
\end{figure}

\par\vspace*{-10mm}\par
\bR
The expressions for $ \kappa _0 $ look complicated. A geometric 
visualization for the simplest cases will be helpful.
Figure 1 depicts the case $ n=3 $. The yellow region is $ {\cal{P}} $, the 
boundary of which is given by the curve $ \kappa =\kappa _0(\lambda ) $. The curve consists of
a straight line (given by (\ref{k01})), of slope $ -n $,
starting from $ -\infty  $, ending at the point $ A $,
and continues along a curvilinear path 
(in red, given by the equation $ \kappa =\lambda +\frac{1}{4\lambda } $). The red curve attains a 
minimum at $ B=\big(\frac{1}{2} ,1\big) $ ; see (\ref{k03}). If we extend the 
rectilinear boundary of $ \cal P $ beyond $ A $, it lies below the red curve, 
a fact manifested in (\ref{k02}).
\eR

\begin{figure}[h]
\begin{center}
\FG{50}{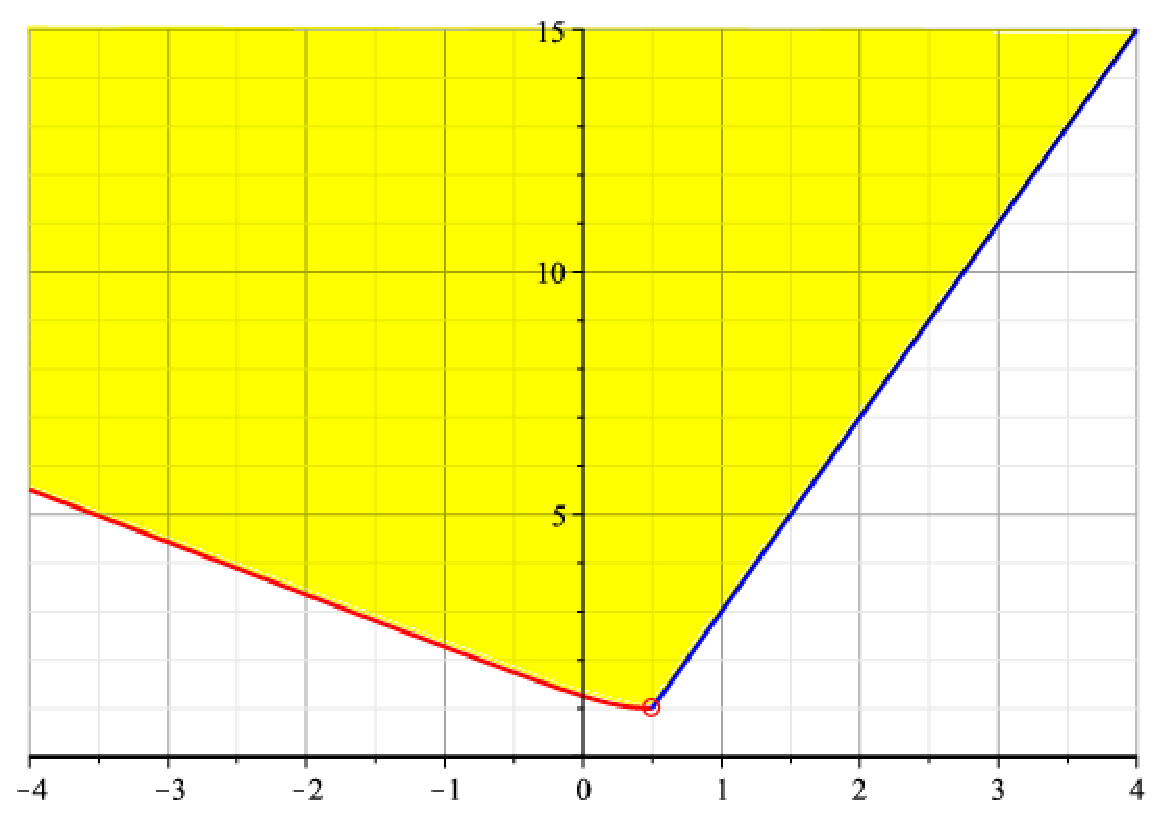}
\end{center}
\par\vspace*{-20mm}\par \hspace*{77mm} {\footnotesize  $ A $}
\par\vspace*{-22mm}\par \hspace*{60mm} $ {\cal{P}}_4 $
\par\vspace*{28mm}\par
\caption{${\cal{P}}_4$, $ A=\big(\frac{1}{2} ,1\big) $}
\par\vspace*{-4mm}\par
\par\vspace*{\baselineskip}\par
\end{figure}

Figure 2 depicts the $ n=4 $ case. 
In contrast to the odd-order case, the rectilinear 
boundary of $ {\cal{P}} $ starts from $ A $ and points up towards
$ \lambda =\infty  $; $ A $ is the lowest point of $ {\cal{P}} $.
The red curvilinear path has a more complicated equation:
$ \textstyle  \kappa  = \frac{9\lambda ^2+9\lambda +2\,\sqrt{(6\lambda ^2-3\lambda +1)^3}-2}{27\lambda ^2} \,. $

In the odd case, the point $ A $ varies according 
to $ n $, while in the even case, $ A $ is fixed. The $ y $-intercept of the 
boundary curve of $ {\cal P} $ is $ \frac{n+1}{2} $ for odd $ n $, but
is always $ \frac{5}{4} $ for even $ n $.

\begin{PROOF}
The existence of a smallest $ \kappa _0 $ for a given $ \lambda  $ follows from the convexity of
$ {\cal{P}} $.

By Lemma 1, 
$ \textstyle  \big(\frac{5}{4} ,  \frac{2n-3}{4n}  \big) \in {\cal P} , $
$ n $ odd. It is well-known that
$$  \theta  = n\sin(x)-\sin(nx)  $$
is NN. Hence, for any $ t>0 $,
$ \phi +t\theta  $ is NN and it corresponds to
$$  \textstyle  \big(nt+\frac{5}{4} ,  \frac{2n-3}{4n} -t \big) \in {\cal P} .  $$
These points are exactly the parametric representation of the straight
line given by (\ref{k01}). We need to show they are actually boundary points.

It is easy to verify that both
$ \phi  $ and $ \theta  $ satisfy (\ref{Bn0}). Hence $ \phi +t\theta  $ satisfies (\ref{Bn0}). Suppose we
decrease the first coefficient $ \kappa  $.
Then (\ref{Bn2}) is no longer satisfied and the polynomial is no longer NN.
Thus, the corresponding $ \kappa  $ must be the smallest possible. 

Similar arguments can be used to prove (\ref{k05})
for the even $ n $ case, by using
$ \sigma  $ of (\ref{s21}) which corresponds to $ \big(1,\frac{1}{2} \big)\in{\cal P} $ (in
place of $ \phi  $) and the NN polynomial
$ n\sin(x)+\sin(nx) $ to construct the rectilinear boundary of $ {\cal P} $.

Next, we look at (\ref{k02}).
Besides satisfying (\ref{Bn0}), $ \phi  $ satisfies (\ref{Bn4}) with equality, which
is equivalent to the easily verified identity
\begin{equation}  1 + 4 \, \sum_{k=1}^{2k} \,(-1)^{k+1} k^3 + (4k-1)(2k+1)^2 = 0 \, .  \Label{k3}  \end{equation}
Let $ \lambda >\frac{2n-3}{4n} $ and $ \kappa =\frac{n+2}{2} -n\lambda  $, then the
polynomial corresponds to
$ \varphi  = \phi  - t\theta  $
for some $ t>0 $. Since both $ \phi  $ and $ \theta  $ satisfy (\ref{Bn0}), so does $ \varphi  $. Now $ \phi  $ satisfies
(\ref{Bn4}) with equality, but $ -t\theta  $ violates (\ref{Bn4}). Hence, $ \varphi  $ violates (\ref{Bn4}). By
Proposition 3, $ \theta  $ cannot be NN.

To prove (\ref{k03}) and (\ref{k04}), we note that direct computation gives
$$  \sigma (x) = \frac{c\big(\frac{1}{2}\big)}{2s\big(\frac{1}{2}\big)} \, (1-c_n).  $$
This implies that $ \sigma \big(\frac{2\pi }{n} \big)=\sigma '\big(\frac{2\pi }{n} \big)=0 $.
For $ \lambda >\frac{1}{2} $,
$$  \eta  = [1,1,\cdots ,\lambda ]_s = \sigma  + \big(\lambda -\textstyle \frac{1}{2}\big)s_n .  $$
It follows that $ \eta \big(\frac{2\pi }{n} \big)=0 $ and $ \eta '\big(\frac{2\pi }{n} \big)=\big(\lambda -\frac{1}{2}\big)c_n\big(\frac{2\pi }{n} \big)<0 $.
Consequently, $ \eta (x) $ is negative in a right neighborhood of $ x=\frac{2\pi }{n} $, and so $ (1,\lambda )\not\in\cal P $.
\end{PROOF}

The equation of the curvilinear boundary of $ {\cal P}_3 $ is determined as follows.
Since
$$  \kappa \sin(x)+\sin(2x)+\lambda \sin(3x) = \sin(x) \left[  4\lambda X^2 +2X +(\kappa -\lambda ) \right]  ,  $$
where $ X=\cos(x) $, NN of the sine polynomial on the LHS
is equivalent to the NN, in $ [-1,1] $, of the algebraic polynomial in $ X $ in square
brackets on the RHS. It is obvious that the latter assertion
is true if $ \kappa  $ is greater than or equal to
$$  - \min_{\scriptscriptstyle X\in[-1,1]} \left\{   4\lambda X^2 +2X -\lambda  \right\}    $$
and the desired conclusion follows.

In theory, the same technique can be used to study the case of general $ n $. The
NN of $ \phi  $ is equivalent to that of an algebraic polynomial
(determined using Chebyshev polynomials). Then $ \kappa _0 $ is the absolute value
of the minimum value of this polynomial in $ [-1,1] $. The determination of
this value, however, becomes more difficult for large $ n $.

The knowledge of $ {\cal P}_3 $ allows us to characterize all NN sine polynomials of
degree 3 and all NN cosine polynomials of degree 2.

\begin{COR}
The sine polynomial $ [a,b,c]_s $ is NN in $ [0,\pi ] $ iff
\begin{itemize}
\item[\rm(i)]  $ |b|\geq 4c $\,\, and \,\, $ a-2|b|+3c\geq 0 $, or
\item[\rm(ii)] $ |b|<4c $\,\, and \,\, $ a\geq c+\frac{b^2}{4c} $.
\end{itemize}
In all cases, a necessary condition is that $ a\geq |b| $.
\end{COR}  \begin{PROOF}
The case $ b=0 $ is trivial.

By making use of the reflection $ x\mapsto-x $, we see that, without loss of generality, we may
assume that $ b>0 $. The general case can be reduced to the case $ b=1 $ by
dividing the sine polynomial by $ b $ and then we are back to the $ {\cal P}_3 $
setting.
\end{PROOF}

\begin{COR}
The cosine polynomial $ a+b\cos(2x)+c\cos(3x) $ is NN in $ [0,\pi ] $ iff
\begin{itemize}
\item[\rm(i)]  $ |b|\geq 4c $\,\, and \,\, $ a-|b|+c\geq 0 $, or
\item[\rm(ii)] $ |b|<4c $\,\, and \,\, $ a\geq c+\frac{b^2}{8c} $.
\end{itemize}
In all cases, a necessary condition is that $ 2a\geq |b|+c $.
\end{COR}  \begin{PROOF}
The identity
$$  a+b\cos(2x)+c\cos(3x) = \frac{(2a-c)\sin(x)+b\sin(2x)+c\sin(3x)}{2\sin(x)}  $$
shows that the NN of the cosine polynomial in question is equivalent to the NN of the sine
polynomial $ [2a-c,b,c]_s $. The conclusions then follow from Corollary 1.
\end{PROOF}

\par\vspace*{-5mm}\par

\end{document}